# SEISMIC INVERSE SCATTERING IN THE 'WAVE-EQUATION' APPROACH

CHRISTIAAN C. STOLK AND MAARTEN V. DE HOOP

ABSTRACT. Seismic data are commonly modeled by a high-frequency single scattering approximation. This amounts to a linearization in the medium coefficient about a smooth background. The discontinuities are contained in the medium perturbation. The wave solutions in the background medium admit a geometrical optics representation. Here we describe the wave propagation in the background medium by a one-way wave equation. Based on this we derive the double-square-root equation, which is a first order pseudodifferential equation, that describes the continuation of seismic data in depth. We consider the modeling operator, its adjoint and reconstruction based on this equation. If the rays in the background that are associated with the reflections due to the perturbation are nowhere horizontal, the singular part of the data is described by the solution to an inhomogeneous double-square-root equation. We derive a microlocal reconstruction equation. The main result is a characterization of the angle transform that generates the common image point gathers, and a proof that this transform contains no artifacts. Finally, pseudodifferential annihilators based on the double-square-root equation are constructed. The double-square-root equation approach is used in seismic data processing.

1. INTRODUCTION

In reflection seismology one places point sources and point receivers on the earth's surface. The source generates acoustic waves in the subsurface, that are reflected where the mediumproperties vary discontinuously. The recorded reflections that can be observed in the data are used to reconstruct these discontinuities.

The data are commonly modeled by a high-frequency single scattering approximation. This amounts to a linearization in the medium coefficient about a smooth background. The discontinuities are contained in the medium perturbation [1]. Thus a linear operator, the modeling operator, depending on the background, that maps the perturbation to the data is obtained. Both the smooth background and the perturbation are in general unknown and have to be reconstructed jointly. In this paper we analyze this reconstruction in the wave-equation approach.

The reconstruction of the perturbation given the background is essentially done by applying the adjoint of mentioned linear map (seismic imaging). The solutions in the background medium admit a geometrical optics representation. Thus the modeling operator is Fourier integral operator [15] (for a general reference of Fourier integral operators see [6]). If the composition of adjoint and modeling operator, the normal operator, is pseudodifferential, then the position of the singularities of the perturbation are recovered by applying

The authors thank The Mathematical Sciences Research Institute for partial support, through NSF grant DMS-9810361. They also thank the members of the institute, and in particular Gunther Uhlmann, for providing a very stimulating environment during the Inverse Problems program in Fall 2001.

the adjoint to the data, and a microlocal reconstruction can be carried out. Under various assumptions on the background, concerning the presence of caustics and the geometry of the rays, results concerning the normal operator have been obtained [1, 8, 12, 10, 16] (section 2.)

In the Kirchhoff approach an approximation to the adjoint operator is constructed that depends on the background only through travel times and complex amplitudes that appear in the geometrical optics approximation. In the so called wave-equation approach [4, 3, 13] the data are downward continued, leading to data from fictitious experiments below the surface at various depths. To form the adjoint a restriction is applied to the downward continued data (imaging condition). Using a one-way wave equation (section 4) for the propagation of waves in the background, the downward continuation is described by the so called double-square-root equation.

The data is formally redundant. There exist invertible Fourier integral operators that generate a set of reconstructions [18] from which the common image point gathers are obtained. Under Beylkin's conditions [1], in particular in the absence of caustics, this can be done by using different subsets of the data. In the presence of caustics the mentioned set can be parameterized by angle between in- and out-going rays at the image point. If the background medium is correct the reconstructions in the set should be the same. This is a criterion that is used in the reconstruction of the background (migration velocity analysis). The redundancy leads to the existence of pseudodifferential operators that annihilate the singular part of the data [18].

In the Kirchhoff approach common image point gathers parameterized by angle can be generated by a generalized Radon transform. In the presence of caustics artifacts were observed in numerical examples in [2]. By microlocal analysis of the Kirchhoff approach the presence of artifacts was shown in [17]. In section 3 an approach to suppress these artifacts is discussed.

In this paper, we consider the modeling, adjoint and reconstruction based on the double-square-root equation approach. The double-square-root equation is a pseudodifferential equation. If the rays in the background that are associated with the reflections due to the perturbation are nowhere horizontal, the singular part of the data is described by the solution to an inhomogeneous double-square-root equation (section 5). We derive a microlocal reconstruction equation in Proposition 6.1. The main result, Theorem 7.1 and Proposition 7.2, is a characterization of the angle transform that generates the common image point gathers, and a proof that this transform contains no artifacts. Annihilators based on the double-square-root equation are constructed in Corollary 8.1.

## 2. High-frequency Born modeling and imaging

We consider the scalar wave equation for acoustic waves in a constant density medium in $\mathbb{R}^n$. In preparation of the later analysis, we distinguish the vertical coordinate $z \in \mathbb{R}$ from the horizontal coordinates $x \in \mathbb{R}^{n-1}$ and write $(x, z) \in \mathbb{R}^n$. In these coordinates the scalar acoustic wave equation is given by

$$(1) \qquad Pu = f, \quad P = c(x,z)^{-2}\frac{\partial^2}{\partial t} + \sum_{j=1}^{n-1} D_{x_j}^2 + D_z^2,$$

where $D_x = -\mathrm{i}\frac{\partial}{\partial x}$, $D_z = -\mathrm{i}\frac{\partial}{\partial z}$. The equation is considered in a time interval $]0, T[$.

If $c \in C^\infty$ the solution operator of (1) propagates singularities along bicharacteristics. These are the solutions of a Hamilton system with Hamiltonian given by the principal symbol of $P$

$$P(x, z, \xi, \zeta, \tau) = -c(x,z)^{-2}\tau^2 + \|\xi\|^2 + \zeta^2.$$

The Hamilton system is given by

$$\text{(2)} \qquad \frac{\partial(x,z,t)}{\partial \lambda} = \frac{\partial P}{\partial(\xi,\zeta,\tau)}, \quad \frac{\partial(\xi,\zeta,\tau)}{\partial \lambda} = -\frac{\partial P}{\partial(x,z,t)}.$$

Its solutions will be parameterized by initial position $(x_0, z_0)$, take-off direction $\alpha \in S^{n-1}$ and frequency $\tau$,

$$\boldsymbol{x} = \boldsymbol{x}(x_0, z_0, \alpha, \tau, t)$$

and similarly for $\boldsymbol{z}, \boldsymbol{t}, \boldsymbol{\xi}, \boldsymbol{\zeta}$; $\tau$ is invariant along the Hamilton flow. Here the evolution parameter is the time $t$. In section 4 we will change the evolution parameter to $z$, and use a similar notation to denote the bicharacteristics parameterized by $z$.

By Duhamel's principle, a causal solution operator for the inhomogeneous equation (1) is given by

$$\text{(3)} \qquad u(x,z,t) = \int_0^t \int G(x,z,t-t_0,x_0,z_0) f(x_0,z_0,t_0)\, \mathrm{d}x_0 \mathrm{d}z_0 \mathrm{d}t_0,$$

where $G$ is a Fourier integral operator with canonical relation that is essentially a union of bicharacteristics. Its kernel can be written as a sum of contributions

$$\text{(4)} \quad G(x,z,t,x_0,z_0) = \sum_i \int_{\mathbb{R}^{N^{(i)}}} a^{(i)}(x,z,t,x_0,z_0,\theta)\, \exp[\mathrm{i}\phi^{(i)}(x,z,x_0,z_0,t,\theta)]\, \mathrm{d}\theta,$$

where the $\phi^{(i)}$ are non-degenerate phase functions and the $a^{(i)}$ suitable symbols, see [6, chapter 5].

We adopt the linearized scattering approximation. The linearization is in the coefficient $c$ around a smooth background $c_0$, $c = c_0 + \delta c$. The perturbation $\delta c$ may contain singularities. We assume that its support is contained in $z > 0$. The perturbation in $G$ at the acquisition surface $z = 0$ is given by

$$\text{(5)} \quad \delta G(r,0,t,s,0) = \int_{\mathbb{R}^{n-1}\times\mathbb{R}_+} \int_0^t G(r,0,t-t_0,x_0,z_0)\, 2c_0^{-3}(x_0,z_0)\delta c(x_0,z_0)$$

$$\partial_{t_0}^2 G(x_0,z_0,t_0,s,0)\, \mathrm{d}t_0 \mathrm{d}x_0 \mathrm{d}z_0,$$

where both $r, s \in \mathbb{R}^{n-1}$. The singular part of $\delta G$ is obtained by substituting (4) into (5). This defines the data modeling map

$$F = F[c_0] \; : \; \delta c \mapsto \mathcal{R}_0 \delta G,$$

where $\mathcal{R}_0$ is the restriction defined by

$$\mathcal{R}_0 \; : \; \mathcal{D}'_\Gamma(\mathbb{R}^{2n+1}) \to \mathcal{D}'(Y)\, , \; u(x,z,t,x_0,z_0) \mapsto (\mathcal{R}_0 u)(x_0,x,t) = u(x,0,t,x_0,0)$$

with acquisition manifold $Y$ a bounded open subset of $\mathbb{R}^{2n-2} \times \mathbb{R}_+$ that contains the range of values of $(s,r,t)$. The restriction is defined on distributions $\mathcal{D}'_\Gamma$ with wavefront sets contained in an open conic set $\Gamma$ such that $\bar{\Gamma}$ does not contain points of the form $(x, 0, t, x_0, 0, 0, \zeta, 0, \zeta_0, 0)$. Since $\tau \neq 0$ on $\mathrm{WF}(\delta G)$ the restriction is well defined. Since $Y$ is bounded and the waves propagate with finite speed we may assume that $\delta c$ is supported in a bounded open subset $X$ of $\mathbb{R}^{n-1} \times \mathbb{R}_+$. We assume that $\overline{X} \cap \{z = 0\} = \emptyset$.

**Assumption 1.** *There are no rays from $(s, 0)$ to $(r, 0)$ with travel time $t$ such that $(s, r, t) \in Y$. For all ray pairs connecting $(r, 0)$ via some $(x, z) \in X$ to $(s, 0)$ with total time $t$ such that $(s, r, t) \in Y$, the rays intersect the plane $z = 0$ transversally at $r$ and $s$.*

**Theorem 2.1.** [15, 10] *With Assumption 1 the map $F$ is a Fourier integral operator $\mathcal{D}'(X) \to \mathcal{D}'(Y)$ of order $(n-1)/4$ with canonical relation*

(6)
$$\big\{(\boldsymbol{x}(x, z, \beta, \tau, t_\mathrm{s}), \boldsymbol{x}(x, z, \alpha, \tau, t_\mathrm{r}), t_\mathrm{s} + t_\mathrm{r}, \boldsymbol{\xi}(x, z, \beta, \tau, t_\mathrm{s}), \boldsymbol{\xi}(x, z, \alpha, \tau, t_\mathrm{r}), \tau; x, z, \xi, \zeta) \mid$$
$$t_\mathrm{s}, t_\mathrm{r} > 0,\ \boldsymbol{z}(x, z, \beta, \tau, t_\mathrm{s}) = \boldsymbol{z}(x, z, \alpha, \tau, t_\mathrm{r}) = 0,\ (\xi, \zeta) = -\tau c_0(x, z)^{-1}(\alpha + \beta),$$
$$(x, z, \alpha, \beta, \tau) \in\ \text{subset of}\ X \times (S^{n-1})^2 \times \mathbb{R}\backslash 0\big\}.$$

Assumption 1 is microlocal. One can identify the set of points $(s, r, t, \sigma, \rho, \tau) \in T^*Y\backslash 0$ where this assumption is violated. If the symbol $\psi = \psi(s, r, t, \sigma, \rho, \tau)$ vanishes on a neighborhood of this set, then the composition $\psi F$ of the pseudodifferential cutoff $\psi = \psi(s, r, t, D_s, D_r, D_t)$ with $F$ is a Fourier integral operator as in the theorem.

We assume that $\psi$ is as before and in addition vanishes outside $Y$. To image the singularities of $\delta c$ from the data we consider the adjoint $F^*\psi$, which is a Fourier integral operator also.

**Assumption 2.** [7] *The projection of the canonical relation (6) on $T^*Y\backslash 0$ is an embedding.*

Since (6) is a canonical relation that projects submersively on the subsurface variables $(x, z, \xi, \zeta)$, the projection of (6) on $T^*Y\backslash 0$ is immersive [9, 25.3.6]. Therefore only the injectivity in the assumption needs to be verified [10].

The following theorem describes the reconstruction of $\delta c$ modulo a pseudodifferential operator with principal symbol that is nonzero at $(x, z, \xi, \zeta)$ whenever there is a point $(s, r, t, \sigma, \rho, \tau; x, z, \xi, \zeta)$ in the canonical relation (6) with $(s, r, t, \sigma, \rho, \tau)$ in the support of $\psi$ (i.e. whenever there is illumination or insonification).

**Theorem 2.2.** *With Assumption 2 the operator $F^*\psi F$ is pseudodifferential of order $n - 1$.*

## 3. GENERALIZED RADON TRANSFORM IN SCATTERING ANGLE

Consider the projection of the canonical relation (6) on the $(x, z, s, r, \tau)$ variables. Where this projection is locally diffeomorphic, the canonical relation (6) can be described by a phase function of the form

$$\tau(T^{(m)}(x, z, s, r) - t)$$

where $T^{(m)}$ is the value of the time variable in (6). There is a set $\{T^{(m)}\}_{m \in M}$ that describes the canonical relation except for a neighborhood of the subset of the canonical relation where mentioned projection is degenerate. Each $T^{(m)}$ is defined on a subset $D^{(m)}$ of $\mathbb{R}^{3n-2}_{(x,z,s,r)}$. We define $F^{(m)}$ to be a contribution to $F$ with phase function given by $T^{(m)}$, such that on a subset of the canonical relation where the projection is nondegenerate $F$ is given microlocally by $\sum_{m \in M} F^{(m)}$.

We can use $(x, z, \xi, \zeta) \in T^*(\mathbb{R}^{n-1} \times \mathbb{R}_+) \backslash 0$ as local coordinates on the canonical relation (6). In addition, we need to parameterize the subsets of the canonical relation given by $(x, z, \xi, \zeta) =$ constant; we denote such parameters by $e$. The canonical relation (6) was parameterized by $(x, z, \alpha, \beta, \tau)$. We relate $(x, z, \xi, \zeta, e)$ by a coordinate transformation to $(x, z, \alpha, \beta, \tau)$: A suitable choice when $\alpha \neq \beta$ is the scattering angles given by

$$(7) \qquad \boldsymbol{e}(x, z, \alpha, \beta) = \left( \arccos(\alpha \cdot \beta), \frac{-\alpha + \beta}{2 \sin(\arccos(\alpha \cdot \beta)/2)} \right) \in ]0, \pi[ \times S^{n-2}.$$

The migration dip $\nu$ is defined as

$$(8) \qquad \nu(\alpha, \beta) = \frac{\alpha + \beta}{\|\alpha + \beta\|} \in S^{n-1}.$$

On $D^{(m)}$ there is a map $(x, z, \alpha, \beta) \mapsto (x, z, s, r)$. We define $\boldsymbol{e}^{(m)} = \boldsymbol{e}^{(m)}(x, z, s, r)$ as the composition of $\boldsymbol{e}$ with the inverse of this map. Likewise, we define $\nu^{(m)} = \nu^{(m)}(x, z, s, r)$.

We define the generalized Radon transform in scattering angle or the Kirchhoff angle transform via a restriction in $F^*$ of the mapping $\boldsymbol{e}^{(m)}$ to a prescribed value $e$, i.e. the distribution kernel of each contribution $F^{(m)*}$ is multiplied by $\delta(e - \boldsymbol{e}^{(m)}(x, z, s, r))$. Its kernel is given by

$$(9) \qquad L(x, z, e, r, s, t) = \sum_{m \in M} (2\pi)^{-(n-1)} \int \overline{A^{(m)}(x, z, s, r, \tau)} e^{i\Phi^{(m)}(x,z,e,s,r,t,\varepsilon,\tau)} \, d\tau d\varepsilon,$$

where $A^{(m)}$ is a symbol for the $m$-th contribution to $F$, supported on $D^{(m)}$, and

$$\Phi^{(m)}(x, z, e, s, r, t, \varepsilon, \tau) = \tau(T^{(m)}(x, z, s, r) - t) + \langle \varepsilon, e - \boldsymbol{e}^{(m)}(x, z, s, r) \rangle.$$

Here, $\varepsilon$ is the cotangent vector corresponding to $e$. Let $\psi_L = \psi_L(D_s, D_r, D_t)$ be a pseudo-differential cutoff such that $\psi(\sigma, \rho, \tau) = 0$ on a conic neighborhood of $\tau = 0$. Then [17] $\psi_L L$ is a Fourier integral operator with canonical relation

$$(10) \quad \cup_{m \in M} \{(x, z, \boldsymbol{e}^{(m)}(x, z, s, r), \boldsymbol{\xi}^{(m)}(x, z, s, r, \tau, \varepsilon), \boldsymbol{\zeta}^{(m)}(x, z, s, r, \tau, \varepsilon), \varepsilon;$$
$$s, r, T^{(m)}(x, z, s, r), \boldsymbol{\rho}^{(m)}(x, z, s, r, \tau, \varepsilon), \boldsymbol{\sigma}^{(m)}(x, z, s, r, \tau, \varepsilon), \tau) \mid$$
$$(x, z, s, r) \in D^{(m)}, \, \varepsilon \in \mathbb{R}^{n-1}, \, \tau \in \mathbb{R} \backslash 0\}$$

where

$$(11) \qquad \boldsymbol{\xi}^{(m)}(x, z, s, r, \tau, \varepsilon) = \partial_x \Phi^{(m)} = \tau \partial_x T^{(m)}(x, z, s, r) - \langle \varepsilon, \partial_x \boldsymbol{e}^{(m)}(x, z, s, r) \rangle,$$

and likewise expressions for $\boldsymbol{\zeta}^{(m)}$, $\boldsymbol{\sigma}^{(m)}$ and $\boldsymbol{\rho}^{(m)}$.

Let $d$ be the Born modeled data. To reveal any artifacts generated by $L$, i.e. singularities in $Ld$ at positions not corresponding to an element of $\mathrm{WF}(\delta c)$, we consider the composition $LF$. This composition is equal to the sum of a smooth $e$-family of pseudodifferential operators and, in general, a non-microlocal operator the wavefront set of which contains no elements with $\varepsilon = 0$ [17, theorem 6.1].

We discuss the practical considerations of the suppression of artifacts generated by $L$ by modification of transform (9). The correct way to remove artifacts from $Ld$ is by a Fourier cut-off in $\varepsilon/\|(\xi,\zeta)\|$ about $0$. For bandlimited data, this is done approximately by local averaging in $e$ at every $(x,z)$. Such procedure has been applied to both synthetic and real data examples in [2]. The artifacts were not fully suppressed. The following modification of the kernel of $L$ (9) was applied,

$$L_\chi(x,z,e,r,s,t) = \sum_{m \in M} (2\pi)^{-(n-1)} \int \overline{A^{(m)}(x,z,s,r,\tau)} \chi(x,z,\nu^{(m)}(x,z,s,r))$$
$$\times e^{i\Phi^{(m)}(x,z,e,s,r,t,\varepsilon,\tau)} \, d\tau d\varepsilon,$$

where $\chi(x,z,\nu)$ is a smooth cutoff function on $\mathbb{R}^n \times S^{n-1}$. Observe that $\nu^{(m)}$ is the direction of $\frac{\partial T^{(m)}}{\partial(x,z)}$.

**Remark 3.1.** The transform $L_\chi$ restricts the wavefront set of operator $L$. If $\mathrm{WF}(\delta c)$ is contained in $\{(x,z,\lambda(\alpha_x(x,z),\alpha_z(x,z))) \mid (x,z,\lambda) \in \mathbb{R}^{n-1} \times \mathbb{R}_+ \times \mathbb{R}\}$, where $(\alpha_x(x,z), \alpha_z(x,z))$ is a smooth covector field on $\mathbb{R}^{n-1} \times \mathbb{R}_+$, there is a $\chi$ with a small support in $\nu$ such that $L_\chi$ generates the true image. Artifacts with singular direction $(\boldsymbol{\xi}^{(m)}, \boldsymbol{\zeta}^{(m)})$ such that the direction of $\frac{\partial T^{(m)}}{\partial(x,z)}$ is outside the support of $\chi$, are suppressed.

The main difficulty arises when the background medium is not (accurately) known. Without knowledge of the background medium there is no criterion to distinguish artifacts from the true image. The approach following the double-square-root equation to be introduced below does not generate artifacts. This is of particular importance in the development of tomographic methods for reconstructing the background medium.

## 4. THE ONE-WAY WAVE OR SINGLE-SQUARE-ROOT EQUATION

In this section we discuss the problem of solving the wave equation by evolution in the vertical, $z$ direction. This problem is in general not well posed, but microlocal solutions can be obtained.

Consider the wave equation rewritten as a first-order system

$$(12) \qquad \frac{\partial}{\partial z}\begin{pmatrix} u \\ \frac{\partial u}{\partial z} \end{pmatrix} = \begin{pmatrix} 0 & 1 \\ -A(x,z,D_x,D_t) & 0 \end{pmatrix} \begin{pmatrix} u \\ \frac{\partial u}{\partial z} \end{pmatrix} + \begin{pmatrix} 0 \\ f \end{pmatrix},$$

where $A(x,z,\xi,\tau) = c_0(x,z)^{-2}\tau^2 - \|\xi\|^2$. Microlocally, away from the zeroes of $A(x,z,\xi,\tau)$, this system can be transformed into diagonal form modulo a smoothing operator [20]. There is a family of pseudodifferential operator matrices $\mathsf{Q}(z) = \mathsf{Q}(x,z,D_x,D_t)$ such that, microlocally,

$$\begin{pmatrix} u_+ \\ u_- \end{pmatrix} = \mathsf{Q}(z)\begin{pmatrix} u \\ \frac{\partial u}{\partial z} \end{pmatrix}, \qquad \begin{pmatrix} f_+ \\ f_- \end{pmatrix} = \mathsf{Q}(z)\begin{pmatrix} 0 \\ f \end{pmatrix},$$

satisfy the one-way wave or single-square-root (SSR) equations

$$(13) \qquad \left(\tfrac{\partial}{\partial z} \pm iB_\pm(x,z,D_x,D_t)\right) u_\pm = f_\pm,$$

see [20]. The principal symbol $b$ of the $B_\pm$ are given by $b(x, z, \xi, \tau) = \sqrt{A(x, z, \xi, \tau)} = \tau\sqrt{\frac{1}{c_0(x,z)^2} - \tau^{-2}\|\xi\|^2}$. For $(x, t, \xi, \tau)$ such that the symbol $B_\pm$ is real, the equation is of hyperbolic type, corresponding (microlocally) to propagating waves. To describe the associated bicharacteristics the Hamiltonian $P$ can be used, as well as $\zeta \mp b(x, z, \xi, \tau)$. The solution of a one-way wave equation describes the propagation of singularities along rays in intervals where $\pm\frac{\partial z}{\partial t} > 0$. The case where the symbol $B_\pm$ is imaginary corresponds to either evanescent waves, or waves that blow up like in a backward heat equation.

We choose the normalization of $Q(z)$, such that (13) is selfadjoint microlocally where the symbol is real,

$$u = Q_+^* u_+ + Q_-^* u_-,$$
$$f_\pm = \pm \tfrac{1}{2}\mathrm{i} Q_\pm f,$$

where $Q_\pm = Q_\pm(z) = Q_\pm(x, z, D_x, D_t)$ are $z$-families of pseudodifferential operators with principal symbols $\tau^{-1/2}\bigl(\frac{1}{c_0(x,z)^2} - \tau^{-2}\|\xi\|^2\bigr)^{-1/4}$.

At the zeroes of $A(x, z, \xi, \tau)$ the operators $B_\pm(z), Q_\pm(z)$ are not yet defined. To this end, we regularize the problem by replacing $A(x, z, \xi, \tau)$ in (12) by

(14) $$c_0(x,z)^{-2}\tau^2 - \|\xi\|^2 - \mathrm{i}\tau^2\phi(x, z, \xi, \tau),$$

where $\phi(x, z, \xi, \tau)$ is positive, small, homogeneous of order 0 in $(\xi, \tau)$ and is supported on a small neighborhood of the set of zeroes of $A(x, z, \xi, \tau)$. With this modification the operators $B_\pm, Q_\pm$ are $z$-families of pseudodifferential operators, defined on the entire cotangent space $T^*\mathbb{R}^n_{(x,t)}\setminus 0$. With this choice of sign for the regularizing imaginary term there is a well defined solution operator $G_-(z, z_0)$, $z_0 > z$, of the initial value problem for $u_-$ given by (13) with $f_- = 0$, see [22, theorem XI.2.1]. The adjoint $G_-(z, z_0)^*$ describes the propagation from $z$ to $z_0$ of (13), regularized in accordance with (14), but with opposite sign of the imaginary part. The operator $G_-$ is a Fourier integral operator with complex phase [11], [9, chapter XXV], [22, chapters X and XI]. By Duhamel's principle a microlocal solution operator for the inhomogeneous equation is given by

(15) $$u_-(\cdot, z) = \int_{-\infty}^z G_-(z, z_0) f_-(\cdot, z_0)\,\mathrm{d}z_0.$$

Microlocally $-\tfrac{1}{2}\mathrm{i}Q_-^*(z) G_-(z, z_0) Q_-(z_0)$ is equal to the Green's function for singularities propagating along bicharacteristics (cf. (2)) with $-\frac{\partial z}{\partial t} > \epsilon$, for some $\epsilon > 0$ that depends on $c_0$ and the support of $\phi$ in (14).

The operator $G_-$ propagates singularities at $(x_0, \xi_0, \tau, z_0)$ along the bicharacteristics for $z$ in an interval $]Z(x_0, \xi_0, \tau, z_0), z_0]$, which is the maximal interval such that the regularized symbol $b$ is real valued. As a consequence the bicharacteristic in this interval is nowhere horizontal. From now we use $z$ as the evolution parameter for bicharacteristics, and denote them as

(16) $(\boldsymbol{x}(x_0, z_0, \xi_0, \tau, z), z, \boldsymbol{t}(x_0, z_0, \xi_0, \tau, z), \boldsymbol{\xi}(x_0, z_0, \xi_0, \tau, z),$
$$b(\boldsymbol{x}(x_0, z_0, \xi_0, \tau, z), z, \boldsymbol{\xi}(x_0, z_0, \xi_0, \tau, z), \tau), \tau).$$

**Remark 4.1.** The symbols of the operators $B_\pm(z), Q_\pm(z)$ can be written as an infinite sum of elementary symbols, $\sum_{i=1}^\infty v_i(x,z,\tau)w_i(z,\xi,\tau)$ say, that is rapidly converging where they are smooth, i.e. away from the set of zeroes of $A(x,z,\xi,\tau)$, in accordance with [21, equation (2.1.11)]. The elementary symbols correspond to multiplications either in horizontal position space ($v_i$) or in horizontal wavenumber (Fourier) space ($w_i$). This is exploited in fast numerical solvers of (13), see [5].

## 5. Modeling, the double-square-root equation

We show that the Born modeling operator can be written, modulo smoothing terms, in terms of the solution operator to the double-square-root (DSR) equation (21) below. We assume that the rays that connect source and receiver to a scattering point in $X$ have nowhere horizontal tangent directions.

**Assumption 3.** *(DSR assumption) If $(x,z) \in X$ and $\alpha, \beta \in S^{n-1}$, $t_s, t_r > 0$ depending on $(x,z,\alpha,\beta)$ are such that $\boldsymbol{z}(x,z,\beta,\tau,t_s) = \boldsymbol{z}(x,z,\alpha,\tau,t_r) = 0$ and $(\boldsymbol{x}(x,z,\beta,\tau,t_s), \boldsymbol{x}(x,z,\alpha,\tau,t_r), t_s + t_r) \in Y$ (cf. (6)), then*

$$\frac{\partial \boldsymbol{z}}{\partial t}(x,z,\beta,\tau,t) < -\epsilon, t \in [0,t_s],$$

$$\frac{\partial \boldsymbol{z}}{\partial t}(x,z,\alpha,\tau,t) < -\epsilon, t \in [0,t_r],$$

*where $\epsilon > 0$ was introduced below (15).*

This assumption is stronger than Assumption 1. This assumption is microlocal, and, given the background medium, a pseudodifferential cutoff can be applied to the data to remove microlocally the part of the data where Assumption 3 is violated.

Under Assumption 3 and the assumption that $\delta c = 0$ on a neighborhood of $z = 0$, the singular part of the Born data is unchanged when $G$ in (5) is replaced by $-\frac{1}{2}\mathrm{i}Q_-^*(x,z,D_x,D_t)G_-Q_-(x,z,D_x,D_t)$. Define the operator $H(z,z_0), z < z_0$ by

$$(17) \quad (H(z,z_0))(s,r,t,s_0,r_0,t_0) =$$
$$\int_\mathbb{R} (G_-(z,z_0))(s,t-t_0-t',s_0,0)(G_-(z,z_0))(r,t',r_0,0)\,\mathrm{d}t'.$$

Here $(G_-(z,z_0))(r,t',r_0,0)$ denotes the distribution kernel of $G_-(z,z_0)$, and similarly for $H(z,z_0)$. Define the maps $\mathcal{I}_1, \mathcal{I}_2$ by

$$(18) \quad \mathcal{I}_1 : \mathcal{D}'(\mathbb{R}^n) \to \mathcal{D}'(\mathbb{R}^{2n-1}) : u(x,z) \mapsto \delta(r-s)u(\tfrac{r+s}{2},z),$$

$$(19) \quad \mathcal{I}_2 : \mathcal{D}'(\mathbb{R}^{2n-1}) \to \mathcal{D}'(\mathbb{R}^{2n}) : u(r,s,z) \mapsto \delta(t)u(\tfrac{r+s}{2},z).$$

The operators $G, G_-, Q_-$ all are of convolution type in the time variable. It follows that the singular part of the Born approximated data (5) is given by

$$(20) \quad F\delta c = \int_0^{z_\mathrm{max}} Q_{-,s}^*(0)Q_{-,r}^*(0)H(0,z)Q_{-,s}(z)Q_{-,r}(z)\tfrac{1}{2}D_t^2(\mathcal{I}_2\mathcal{I}_1\,c_0^{-3}\delta c)(\cdot,z)\,\mathrm{d}z.$$

where $Q_{-,s}(z)$ is short for $Q_-(s,z,D_s,D_t)$ and similarly for $Q_{-,r}(z)$. Here $z_\mathrm{max}$ is the maximum depth illuminated from acquisition manifold $Y$ given the background medium $c_0$.

Define the inhomogeneous double-square-root (DSR) equation by

$$(21) \quad (\tfrac{\partial}{\partial z} - \mathrm{i}B_-(s,z,D_s,D_t) - \mathrm{i}B_-(r,z,D_r,D_t))u = g(s,r,t,z).$$

It follows from the definition of $G_-(z,z_0)$ that the operator $H(z,z_0)$ is a solution operator for the Cauchy initial value problem for the (regularized) DSR equation. The Born approximated data is modulo smooth functions equal to the solution of the DSR equation at $z = 0$ with inhomogeneous term $g(s,r,t,z) = \tfrac{1}{2}D_t^2(\mathcal{I}_2\mathcal{I}_1\, c_0^{-3}\delta c)(s,r,t,z)$.

The operator $H(z,z_0)$ propagates singularities along the bicharacteristics for $z$ in an interval $]Z_{\min}(s_0,r_0,\sigma_0,\rho_0,\tau),z_0]$, which is the intersection of the one-way intervals associated with the source and receiver bicharacteristics, in the notation of (16),

$$(22) \quad (\boldsymbol{x}(s_0,z_0,\sigma_0,\tau,z), \boldsymbol{x}(r_0,z_0,\rho_0,\tau,z), t_0 + \boldsymbol{t}(s_0,z_0,\sigma_0,\tau,z) + \boldsymbol{t}(r_0,z_0,\rho_0,\tau,z), z,$$

$$\boldsymbol{\xi}(s_0,z_0,\sigma_0,\tau,z), \boldsymbol{\xi}(r_0,z_0,\rho_0,\tau,z), \tau,$$

$$\boldsymbol{\Gamma}(\boldsymbol{x}(s_0,z_0,\sigma_0,\tau,z), \boldsymbol{x}(r_0,z_0,\rho_0,\tau,z), \boldsymbol{\xi}(s_0,z_0,\sigma_0,\tau,z), \boldsymbol{\xi}(r_0,z_0,\rho_0,\tau,z), \tau, z)),$$

where

$$(23) \quad \boldsymbol{\Gamma}(s,r,\sigma,\rho,\tau,z) = b(s,z,\sigma,\tau) + b(r,z,\rho,\tau).$$

## 6. Double-square-root reconstruction

Let $N := F^*\psi F$, where $\psi = \psi(s,r,t,D_s,D_r,D_t)$ is a suitable pseudodifferential cutoff as in section 2. By Theorem 2.1, with Assumption 2, $N = N(x,z,D_x,D_z)$ is a pseudodifferential operator. Therefore

$$(24) \quad N(x,z,D_x,D_z)\delta c = F^*\psi(s,r,t,D_s,D_r,D_t)d,$$

and $\delta c$ can be reconstructed microlocally where the principal symbol of $N$ is non-zero. Assumption 3 implies the injectivity in Assumption 2 and is hence stronger than Assumption 2 by the remark below Assumption 2. We present a reconstruction formula similar to (24) based on the DSR Born modeling (20). This leads, again, to reconstruction modulo a pseudodifferential operator for which an explicit expression is given.

The adjoint operator $H(0,z)^*$ propagates the data downward and backward in time. We consider $H(0,z)^*d$ as a function of $(s,r,t,z)$. The operator $H(0,z)^*$ is microlocally a Fourier integral operator with real phase and canonical relation which follows from (22) with $z_0 = 0$, for $z$ in an interval $[0, Z_{\max}(s_0,r_0,\sigma_0,\rho_0,\tau)[$

$$(25) \quad \{(\boldsymbol{x}(s_0,0,\sigma_0,\tau,z), \boldsymbol{x}(r_0,0,\rho_0,\tau,z), t_0 + \boldsymbol{t}(s_0,0,\sigma_0,\tau,z) + \boldsymbol{t}(r_0,0,\rho_0,\tau,z), z,$$

$$\boldsymbol{\xi}(s_0,0,\sigma_0,\tau,z), \boldsymbol{\xi}(r_0,0,\rho_0,\tau,z), \tau, \boldsymbol{\Gamma}(s_0,r_0,\sigma_0,\rho_0,\tau,z); s_0,r_0,t_0,\sigma_0,\rho_0,\tau) \mid$$

$$(s_0,r_0,t_0,\sigma_0,\rho_0,\tau) \in T^*\mathbb{R}^{2n-1}\backslash 0, z \in [0, Z_{\max}(s_0,r_0,\sigma_0,\rho_0,\tau)[\}.$$

The adjoint of the operator $\mathcal{I}_2$ (cf. (19)) is given by the restriction $\mathcal{R}_2$ defined by

$$g(r,s,t,z) \mapsto (\mathcal{R}_2 g)(r,s,z) = g(r,s,0,z).$$

If $u = u(s,r,z)$ then $Ku = (Ku)(s,r,t)$ will be defined by

$$Ku = \int_{-\infty}^{0} H(0,z)\mathcal{I}_2 u(\cdot,z)\,\mathrm{d}z.$$

Let $\psi$ be a pseudodifferential cutoff in the $(s, r, t)$ variables supported in

$$\{(s_0, r_0, t_0, \sigma_0, \rho_0, \tau) \,|\, t_0 \in$$
$$[0, -\boldsymbol{t}(s_0, 0, \sigma_0, \tau, Z_{\max}(s_0, r_0, \sigma_0, \rho_0, \tau)) - \boldsymbol{t}(r_0, 0, \rho_0, \tau, Z_{\max}(s_0, r_0, \sigma_0, \rho_0, \tau))[\}.$$

Then $\psi K$ is a Fourier integral operator with real phase. The adjoint $K^*$ of $K$ follows from the equality

$$\langle \psi d, \int_{-\infty}^0 H(0,z) \mathcal{I}_2 u(\cdot, z) \,\mathrm{d}z \rangle_{(s,r,t)} = \int_{-\infty}^0 \langle \mathcal{R}_2 H(0,z)^* \psi d, u(\cdot, z) \rangle_{(s,r)} \,\mathrm{d}z$$

and is given by

$$K^* \psi = \mathcal{R}_2 H(0,z)^* \psi.$$

Since $\frac{\partial \boldsymbol{t}}{\partial z} < 0$, $\boldsymbol{t}$ depends strictly monotone on $z$ for each $(s_0, r_0, t_0, \sigma_0, \rho_0, \tau)$. The DSR bicharacteristic with initial values $(s_0, r_0, t_0, \sigma_0, \rho_0, \tau)$ hence intersects the $t = 0$ hyperplane at most once, and the intersection is transversal. From this, it follows that the composition $K^* \psi$ is a Fourier integral operators with a locally invertible canonical relation.

The canonical relation maps points $(s_0, r_0, t_0, \sigma_0, \rho_0, \tau)$ in a subset of the cotangent acquisition space $T^* \mathbb{R}^{2n-1} \backslash 0$ to $(r, s, z, \sigma, \rho, \zeta)$ in a subset of $T^* \mathbb{R}^{2n-1} \backslash 0$. This map converts time to depth. We define the symbol $\Psi(s, r, z, \sigma, \rho, \zeta)$ as the pull back of $\psi(s_0, r_0, t_0, \sigma_0, \rho_0, \tau)$ under the inverse of this map. Starting from the Born modeling (20), we find the following reconstruction proposition.

**Proposition 6.1.** *There are pseudodifferential operators $\Phi = \Phi(x, z, D_x, D_z)$ of order $n - 1$ with principal symbol*

$$(26) \qquad \Phi(x, z, \xi, \zeta) = \int_{\mathbb{R}^{n-1}} \Psi(x, x, z, \tfrac{1}{2}\xi - \theta, \tfrac{1}{2}\xi + \theta, \zeta) \,\mathrm{d}\theta,$$

*and $\Xi(z) = \Xi(s, r, t, D_s, D_r, D_t, z)$ of order $0$ with principal symbol*

$$\Xi(s, r, t, \sigma, \rho, \tau, z) = \left| \frac{\partial \boldsymbol{\Gamma}}{\partial \tau}(s, r, \sigma, \rho, \tau, z) \right|$$
$$= c_0(s,z)^{-2}(c_0(s,z)^{-2} - \tau^{-2}\|\sigma\|^2)^{-1/2} + c_0(r,z)^{-2}(c_0(r,z)^{-2} - \tau^{-2}\|\rho\|^2)^{-1/2},$$

*such that*

$$(27) \quad \Phi(x, z, D_x, D_z) \delta c$$
$$= 2c_0^3 \mathcal{R}_1 \mathcal{R}_2 \Xi(z) Q_{-,s}^*(z)^{-1} Q_{-,r}^*(z)^{-1} H(0,z)^* Q_{-,s}(0)^{-1} Q_{-,r}(0)^{-1} D_t^{-2} \psi d,$$

*where $d = F \delta c$ is the Born modeled data.*

*Proof.* We calculate microlocally the principal symbol of $K^* K$. The kernel of the operator $H(0, z)$ has microlocally an oscillatory integral representation with a phase function associated with generating function, $S = S(z, s, r, t, y_{0I}, \eta_{0J})$ say, where $y_0 = (s_0, r_0, t_0)$ and

$\eta_0$ is the corresponding cotangent vector, and $\{I, J\}$ is a partition of $\{1, \ldots, 2n-1\}$,

$$H(0,z)(y_0, s, r, t) = (2\pi)^{-(2n-1+|I|)/2}$$
$$\times \int A(z, s, r, t, y_0, \eta_{0J}) e^{i(S(z,s,r,t,y_{0I},\eta_{0J}) - \langle \eta_{0J}, y_{0J} \rangle)} \, d\eta_{0J},$$

The adjoint $H(0,z)^*$ has amplitude $\overline{A(z, s, r, t, y_0, \eta_{0J})}$ and phase $-S(z, s, r, t, y_{0I}, \eta_{0J}) + \langle \eta_{0J}, y_{0J} \rangle$. Hence, the kernel of the composition $H(0,z)^* H(0,z)$ has the oscillatory integral representation

$$(2\pi)^{-(2n-1)} \int \overline{A(z, s', r', t', y_0, \eta_{0J})} A(z, s, r, t, y_0, \eta_{0J})$$
$$\times e^{i[-S(z,s',r',t',y_{0I},\eta_{0J}) + S(z,s,r,t,y_{0I},\eta_{0J})]} \, dy_{0I} d\eta_{0J}.$$

We expand the phase in a Taylor series around $(s', r', t') = (s, r, t)$ and identify the gradient

$$-\frac{\partial S}{\partial(s, r, t)}(z, s, r, t, y_{0I}, \eta_{0J}) =$$
$$(\sigma(z, s, r, t, y_{0I}, \eta_{0J}), \rho(z, s, r, t, y_{0I}, \eta_{0J}), \tau(z, s, r, t, y_{0I}, \eta_{0J})).$$

Applying a change of variables, $(y_{0I}, \eta_{0J}) \mapsto (\sigma, \rho, \tau)$, the phase takes the form

$$\langle (\sigma, \rho, \tau), (s' - s, r' - r, t' - t) \rangle.$$

Since $H(0,z)^* H(0,z)$ is a pseudodifferential operator with symbol 1, microlocally in the support of the cutoff $\psi$, we conclude that the principal part $a$ of the amplitude $A$ satisfies

$$|a(z, s, r, t, y_0, \eta_{0J})| = \left| \frac{\partial(\sigma, \rho, \tau)}{\partial(y_{0I}, \eta_{0J})} \right|^{1/2}.$$

The kernel of operator $K$ has an oscillatory integral representation similar to the one of $H(0,z)$,

$$K(y_0, s, r, z) = (2\pi)^{-(2n-1+|I|)/2} \int A(z, s, r, 0, y_0, \eta_{0J}) e^{i(S(z,s,r,0,y_{0I},\eta_{0J}) - \langle \eta_{0J}, y_{0J} \rangle)} \, d\eta_{0J}$$

(we have applied $\mathcal{I}_2$ at $z$). It follows that the composition $K^*K$, carrying out an analysis similar to the one for $H(0,z)^* H(0,z)$, is a pseudodifferential operator, microlocally. Its amplitude has principal part

$$\left| \frac{\partial(\zeta, \sigma, \rho)}{\partial(y_{0I}, \eta_{0J})} \right|^{-1} \left| \frac{\partial(\sigma, \rho, \tau)}{\partial(y_{0I}, \eta_{0J})} \right|.$$

For fixed $(s, r, \sigma, \rho, z)$ the map $\tau \mapsto \mathbf{\Gamma}(s, r, \sigma, \rho, \tau, z)$ is invertible on a set given by $|\tau|$ sufficiently large. This map will be denoted by $\mathbf{\Gamma}^{-1} = \mathbf{\Gamma}^{-1}(s, r, z, \sigma, \rho, \zeta)$. It follows that the principal part of $K^*K$ is microlocally given by

$$\left| \frac{\partial \mathbf{\Gamma}}{\partial \tau}(s, r, \sigma, \rho, \mathbf{\Gamma}^{-1}(s, r, z, \sigma, \rho, \zeta), z) \right|^{-1}.$$

We finally consider the composition of operators $\mathcal{R}_1 \Psi \mathcal{I}_1$. Its kernel has an oscillatory integral representation,

$$(2\pi)^{-(2n-1)} \int_{\mathbb{R}^{2n-1}} \Psi(x, x, z, \sigma, \rho, \zeta) \, e^{i\langle (x,x,z)-(x',x',z'),(\sigma,\rho,\zeta)\rangle} \, d\rho d\sigma d\zeta =$$

$$(2\pi)^{-(2n-1)} \int_{\mathbb{R}^n} \int_{\mathbb{R}^{n-1}} \Psi(x, x, z, \tfrac{1}{2}\xi - \theta, \tfrac{1}{2}\xi + \theta, \zeta) \, d\theta \, e^{i(\langle (x-x'),\xi\rangle + \langle (z-z'),\zeta\rangle)} \, d\xi d\zeta,$$

upon changing variables of integration, $\sigma = \tfrac{1}{2}\xi - \theta$, $\rho = \tfrac{1}{2}\xi + \theta$. The domain of the $\theta$ integral is bounded depending on $(\xi, \zeta)$, since $\Psi$ is a cutoff in $(s, r, t, \sigma, \rho, \tau)$. Hence $\mathcal{R}_1 \Psi \mathcal{I}_1$ is a pseudodifferential operator of order $n - 1$ with principal symbol (26).

Applying the above results to the expression (20) for the Born modeled data leads to the statement of the proposition. $\square$

**Remark 6.2.** Note that in (27) all the operators on the right hand side, except the downward continuation $H$, act at depth $z$ only. On the contrary, the operator $\Phi(x, z, D_x, D_z)$ on the left hand side depends on the Hamiltonian flow associated with the background medium in the depth interval $[0, z]$. In the usual wave-equation imaging algorithms it is hence straightforward to include the pseudodifferential factors on the right hand side of (27). To account for the operator $\Phi(x, z, D_x, D_z)$ on the left hand side one requires an additional ray computation.

**Remark 6.3.** Depending on the background medium, the reconstruction can also be done using data on a submanifold $Y'$ of $Y$. Let $\mathcal{R}'$ be the restriction of a function on $Y$ to $Y'$, so that the forward map for this case is given by $\mathcal{R}'F$. In suitable local coordinates $(y', y'')$ on $Y$ such that $y'' = 0$ defines $Y'$, the adjoint $\mathcal{I}'$ of $\mathcal{R}'$ is given by the map $(\mathcal{I}'f)(y', y'') = f(y')\delta(y'')$. Conditions such that $F^*\mathcal{I}'\psi'\mathcal{R}'F$ is pseudodifferential are given in [12], where $\psi'$ is a suitable pseudodifferential cutoff. Reconstruction modulo a pseudodifferential operator is done in this case by first applying the map $\mathcal{I}'$ to the data, and then applying the previous procedure. Applying $\mathcal{I}'$ to the data simply means adding zeroes where there is no data in $Y$.

## 7. THE WAVE-EQUATION ANGLE TRANSFORM

We define the wave-equation angle transform $A_{\text{WE}}$ by the following integral of the downward continued data, $H(0, z)^*d$,

(28) $\quad (A_{\text{WE}}d)(x, z, p) = \int_{\mathbb{R}^{n-1}} (H(0, z)^*\psi d)(x - \tfrac{h}{2}, x + \tfrac{h}{2}, ph)\chi(x, z, h) \, dh,$

(cf. [14]), where $h \mapsto \chi(x, z, h)$ is a compactly supported cutoff function the support of which contains $h = 0$.

**Theorem 7.1.** *Let $C_0$ be an upper bound for $c_0$. Assume that*

(29) $\quad\quad\quad\quad\quad\quad\quad\quad |p| < p_{\max} < \tfrac{1}{2}C_0^{-1}.$

*Then $A_{\text{WE}}$ is a Fourier integral operator such that $A_{\text{WE}}F$ is a smooth $p$-family of pseudodifferential operators in $(x, z)$. Let $C_1$ be an upper bound for $\frac{\partial c_0^{-2}}{\partial x}$. If in addition the* *function $h \mapsto \chi(x, z, h)$ is supported in $B(0, R)$, where $R$ depends on $C_1$ and $\psi$, then the canonical relation of $A_{\text{WE}}$ corresponds to an invertible map from a subset of $T^*\mathbb{R}^{2n-1}_{(s,r,t)}$ to a subset of $T^*\mathbb{R}^{2n-1}_{(x,z,p)}$ that has nonempty intersection with the set $\theta = 0$ (where $\theta$ is the p-covector).*

*Proof.* The map $d \mapsto H(0, z)^*\psi d$ is a Fourier integral operator with canonical relation given in (25).

The Schwarz kernel of the map $H(0, z)^*\psi d \mapsto A_{\text{WE}} d$ equals

$$\delta(x - \tfrac{s+r}{2})\, \delta(p(r-s) - t)\, \delta(z - z')\, \chi(x, z, r-s)$$
$$= (2\pi)^{-n-1} \int e^{i(\langle \xi, x - \frac{s+r}{2}\rangle + \tau(p(r-s)-t) + \zeta(z-z'))}\, d\xi\, d\tau\, d\zeta.$$

It is a Fourier integral operator with canonical relation that is contained in $T^*\mathbb{R}^{2n-1}_{(x,z,p)} \setminus 0 \times T^*\mathbb{R}^{2n}_{(s,r,t,z)} \setminus 0$ and given by

(30)   $\{(\tfrac{s+r}{2}, z, p, \xi, \zeta, (r-s)\tau; s, r, p(r-s), z, \tfrac{\xi}{2} + p\tau, \tfrac{\xi}{2} - p\tau, \tau, \zeta) \,|$
$$(s, r, z, p, \xi, \zeta, \tau) \in \text{(subset of )}\mathbb{R}^{4n-1}\}.$$

This canonical relation can be parameterized by the coordinates of $T^*\mathbb{R}^{2n}_{(s,r,t,z)} \setminus 0$ except $t$, that is $(s, r, z, \sigma, \rho, \tau, \zeta)$. The projection of (30) on $T^*\mathbb{R}^{2n}_{(s,r,t,z)} \setminus 0$ is a hypersurface defined by

(31) $$t = \left\langle \frac{\sigma - \rho}{2\tau}, (r-s) \right\rangle.$$

The canonical relation of the map $d \mapsto H(0, z)^*\psi d$, considered as a function of $(s, r, t, z)$, is parameterized by $(s_0, r_0, t_0, \sigma_0, \rho_0, \tau, z)$. The canonical relation (25) is time translation invariant and the line in $T^*\mathbb{R}^{2n}_{(s,r,t,z)} \setminus 0$ parameterized by $t_0$ for fixed $(s_0, r_0, \sigma_0, \rho_0, \tau)$ intersects the hypersurface (31) transversally. It follows that the composition of the canonical relations (25) and (30) is transversal. The composition is parameterized by $(s_0, r_0, \sigma_0, \rho_0, \tau, z)$. It follows that $A_{\text{WE}}$ is a Fourier integral operator.

The composition $H(0, z)^*\psi F$, that maps $\delta c = \delta c(x, z')$ to the downward continued data as a function of $(s, r, t, z)$, is a Fourier integral operator with canonical relation

(32)   $\{(\boldsymbol{x}(x, z', \sigma', \tau, z), \boldsymbol{x}(x, z', \rho', \tau, z), \boldsymbol{t}(x, z', \sigma', \tau, z) + \boldsymbol{t}(x, z', \rho', \tau, z), z,$
$$\boldsymbol{\xi}(x, z', \sigma', \tau, z), \boldsymbol{\xi}(x, z', \rho', \tau, z), \tau, \zeta; x, z', \sigma' + \rho', \zeta') \,|$$
$$(x, z', \sigma', \rho', \zeta', z) \in \text{ a subset of } \mathbb{R}^{3n}, \tau \text{ such that } \zeta' = \boldsymbol{\Gamma}(x, x, \sigma', \rho', \tau, z'),$$
$$\zeta = \boldsymbol{\Gamma}(\boldsymbol{x}(x, z', \sigma', \tau, z), \boldsymbol{x}(x, z', \rho', \tau, z), \boldsymbol{\xi}(x, z', \sigma', \tau, z), \boldsymbol{\xi}(x, z', \rho', \tau, z), \tau, z)\}.$$

The propagation of singularities upward by $F$ and downward by $H(0, z)^*$ is along the same DSR bicharacteristics.

We show that the composition $A_{\text{WE}} F$ is a Fourier integral with canonical relation contained in

(33)   $\{(x, z, p, \xi, \zeta, 0; x, z, \xi, \zeta) \,|\, (x, z, \xi, \zeta) \in T^*\mathbb{R}^n_{(x,z)} \setminus 0, p \in\, ]-p_{\max}, p_{\max}[\}.$

From that it follows that $A_{\text{WE}}F$ is a $p$-family of pseudodifferential operators. The projection of (32) on $T^*\mathbb{R}^{2n}_{(s,r,t,z)}$ intersects the hypersurface (31) at $z' = z$, since then $r - s = 0$ and $t = 0$, leading to elements in (33). Since singularities propagate with speed less than $C_0$, if $(v_s, v_r, v_t, v_z, v_\sigma, v_\rho, 0, v_\zeta)$ is a tangent vector to the DSR bicharacteristic, then $\frac{v_s - v_r}{v_t} \leq 2C_0$. Therefore, by (29), the composition of (32) with (30) is transversal and contains no elements outside (33).

The projection of the canonical relation of $A_{\text{WE}}$ on the second component $T^*\mathbb{R}^{2n-1}_{(s,r,t)}\backslash 0$ is invertible if each DSR bicharacteristic with initial values $(s_0, r_0, t_0, \sigma_0, \rho_0, \tau)$, parameterized by $z$, intersects the hypersurface (31) at most once and transversally. Let $\boldsymbol{p}(z)$ denote $\frac{\sigma - \rho}{2\tau}$ along a certain DSR bicharacteristic and let $\boldsymbol{t}(z)$ denote the time and $\boldsymbol{h}(z)$ denote the value of $r - s$. The elements of the canonical relation of $A_{\text{WE}}$ correspond to solutions of $\boldsymbol{t}(z) - \langle \boldsymbol{p}(z), \boldsymbol{h}(z)\rangle = 0$. To estimate the derivative of the left hand side we observe that

$$\frac{\partial \boldsymbol{t}}{\partial z} - \langle \boldsymbol{p}(z), \frac{\partial \boldsymbol{h}}{\partial z}(z)\rangle < -\epsilon_0$$

for some $\epsilon_0 > 0$ depending on the cutoff $\psi$ (or $\phi$ in (14)) and on the value $p_{\max}$. Since

$$\frac{\partial \boldsymbol{\xi}}{\partial z} = \frac{\partial b}{\partial x} = -\frac{\tau}{\sqrt{c_0^{-2} - \tau^{-2}\|\xi\|^2}} \frac{\partial c_0^{-2}}{\partial x},$$

there is a bound on $\frac{\partial \boldsymbol{p}}{\partial z}$ in terms of $C_1$ and on $\psi$ (to bound the square root from below). It follows that for some $\epsilon_1 < \epsilon_0$

$$\left|\langle \frac{\partial \boldsymbol{p}}{\partial z}(z), \boldsymbol{h}(z)\rangle\right| < \epsilon_1 < \epsilon_0$$

if $\|h\| < C_2 C_1^{-1}$ for some constant $C_2$ depending on $\psi$ and $p_{\max}$. This implies that the function $z \mapsto \boldsymbol{t}(z) - \langle \boldsymbol{p}(z), \boldsymbol{h}(z)\rangle$ is monotone. Hence the projection of the canonical relation of $A_{\text{WE}}$ on $T^*\mathbb{R}^{2n-1}_{(s,r,t)}\backslash 0$ is invertible. It follows from the above reasoning that the projection on $T^*\mathbb{R}^{2n-1}_{(x,z,p)}\backslash 0$ is invertible as well. This establishes the last statement of the theorem. □

To conclude this section we determine, at the principal symbol level, the modification of (28) that leads to microlocal reconstruction.

**Proposition 7.2.** *Define $\widetilde{A}_{\text{WE}}$ by*

$$(\widetilde{A}_{\text{WE}}d)(x, z, p) = \int_{\mathbb{R}^{n-1}} \chi(x, z, h) 2c_0(x, z)^3$$
$$\times \left(\Xi(z) Q_{-,s}^*(z)^{-1} Q_{-,r}^*(z)^{-1} H(0, z)^* Q_{-,s}(0)^{-1} Q_{-,r}(0)^{-1} D_t^{-2} \psi d\right)\left(x - \tfrac{h}{2}, x + \tfrac{h}{2}, ph\right) \mathrm{d}h.$$

*Suppose that $\chi(x, z, 0) = 1$ and $h \mapsto \chi(x, z, h)$ is supported in $B(0, R)$ (cf. Theorem 7.1), then $\widetilde{A}_{\text{WE}}$ is an invertible Fourier integral operator. Let the symbol $\Psi_{\text{WE}} = \Psi_{\text{WE}}(x, z, p, \xi, \zeta, \theta)$ (where $\theta$ is the $p$-covector) be given by the pull back of $\psi$ under the map from $T^*\mathbb{R}^{2n-1}_{(x,z,p)}\backslash 0$ to $T^*\mathbb{R}^{2n-1}_{(r,s,t)}\backslash 0$ induced by the canonical relation of $A_{\text{WE}}$. The composition $\widetilde{A}_{\text{WE}}F$ is a $p$-family of pseudodifferential operators with principal symbol $\Psi_{\text{WE}}(x, z, p, \xi, \zeta, 0)$.*

*Proof.* It is sufficient to show that the map

$$\delta c \mapsto \int \left( \Xi(z) H(0,z)^* \int H(0,z') \mathcal{I}_2 \mathcal{I}_1 \delta c \, \mathrm{d}z' \right) (x - \tfrac{h}{2}, x + \tfrac{h}{2}, ph) \, \mathrm{d}h, \tag{34}$$

microlocally has principal symbol equal to 1. In this proof we will omit the cutoff functions that are part of the symbols; the calculations will be valid microlocally on the support of a cutoff.

Using an oscillatory integral representation of $H$ similar to the one in the proof of Proposition 6.1, we find that the principal contribution to the kernel of this map, as a function of $(x, z, p; x', z')$, can be written as

$$(2\pi)^{-(2n-1)} \int \Xi(x - \tfrac{1}{2}h, x + \tfrac{1}{2}h, ph, -\frac{\partial S}{\partial s}, -\frac{\partial S}{\partial r}, -\frac{\partial S}{\partial t}, z)$$
$$\times \overline{A(z, x - \tfrac{1}{2}h, x + \tfrac{1}{2}h, ph, y_0, \eta_{0J})} A(z', x', x', 0, y_0, \eta_{0J})$$
$$\times \mathrm{e}^{\mathrm{i}[-S(z, x - \tfrac{1}{2}h, x + \tfrac{1}{2}h, ph, y_{0I}, \eta_{0J}) + S(z', x', x', 0, y_{0I}, \eta_{0J})]} \, \mathrm{d}y_{0I} \mathrm{d}\eta_{0J} \, \mathrm{d}h.$$

We expand the phase in a Taylor series around $(x', z', h) = (x, z, 0)$ and identify the gradient at $(x, z, 0)$,

$$-\frac{\partial S}{\partial x}(z, x, x, 0, y_{0I}, \eta_{0J}) = \sigma(z, x, x, 0, y_{0I}, \eta_{0J}) + \rho(z, x, x, 0, y_{0I}, \eta_{0J}),$$
$$-\frac{\partial S}{\partial z}(z, x, x, 0, y_{0I}, \eta_{0J}) = \zeta(z, x, x, 0, y_{0I}, \eta_{0J}),$$
$$-\frac{\partial S}{\partial h}(z, x, x, 0, y_{0I}, \eta_{0J})$$
$$= -\tfrac{1}{2}\sigma(z, x, x, 0, y_{0I}, \eta_{0J}) + \tfrac{1}{2}\rho(z, x, x, 0, y_{0I}, \eta_{0J}) + p\tau(z, x, x, 0, y_{0I}, \eta_{0J}).$$

Applying a change of variables, $(y_{0I}, \eta_{0J}) \mapsto (\sigma, \rho, \zeta)$, the phase takes the form

$$\langle \sigma + \rho, x - x' \rangle + \zeta(z - z') + \langle \tfrac{1}{2}(\rho - \sigma), h \rangle.$$

The amplitude factor $\Xi \overline{A} A$ becomes equal to one by the calculations in the proof of Proposition 6.1. Upon changing integration variables, $\sigma = \tfrac{1}{2}\xi - \theta, \rho = \tfrac{1}{2}\xi + \theta$, the oscillatory integral takes the leading-order form (microlocally)

$$(2\pi)^{-(2n-1)} \int \mathrm{e}^{\mathrm{i}(\langle \xi, x-x' \rangle + \zeta(z-z') + \langle \theta, h \rangle)} \, \mathrm{d}\xi \mathrm{d}\zeta \mathrm{d}\theta \mathrm{d}h.$$

It follows by integrating the $\theta, h$ variables that (34) indeed has principal symbol equal to 1 microlocally. □

## 8. ANNIHILATORS

It was observed in [18] that there are pseudodifferential operators that annihilate the data, due to the fact that the inverse problem is formally overdetermined. On transformed data $\widetilde{A}_{\mathrm{WE}} d$, these are given by $\frac{\partial}{\partial p_i}, i = 1, \ldots, n-1$ hence the annihilators are given by

$$\langle \widetilde{A}_{\mathrm{WE}} \rangle^{-1} \frac{\partial}{\partial p_i} \widetilde{A}_{\mathrm{WE}},$$

Where $\langle \widetilde{A}_{\text{WE}} \rangle^{-1}$ is a regularized inverse for $\widetilde{A}_{\text{WE}}$, that is a microlocal inverse for a subset of $T^*\mathbb{R}^{2n-1}_{(x,z,p)}$ where $\widetilde{A}_{\text{WE}}$ is invertible.

We define the operator $\widetilde{K}^*$ that maps $d$ to a function $(s, r, z)$ by

$$(\widetilde{K}^*d)(s,r,z) = 2c_0((s+r)/2, z)^3$$
$$\times \left( \mathcal{R}_2 \Xi(z) Q^*_{-,s}(z)^{-1} Q^*_{-,r}(z)^{-1} H(0,z)^* Q_{-,s}(0)^{-1} Q_{-,r}(0)^{-1} D_t^{-2} d \right)(s,r,z).$$

Using (27), we observe that the operator $\widetilde{K}^*$ acting on $d = F\delta c$ yields

$$\widetilde{K}^*\psi d = \Psi(s, r, z, D_s, D_r, D_z) \mathcal{I}_1 \delta c$$

(for the definition of $\Psi$ see remark above Proposition 6.1). Applying the operator $\mathcal{M}$ given by the multiplication by $r - s$ to this equation yields

$$\mathcal{M}\widetilde{K}^*\psi d = [\mathcal{M}, \Psi(s, r, z, D_s, D_r, D_z)] \mathcal{I}_1 \delta c,$$

i.e. there is only a lower order contribution, since $(r-s)\delta(r-s) = 0$ and hence $\mathcal{M}\mathcal{I}_1 = 0$ (cf. (18)).

For the subset in $T^*\mathbb{R}^{2n-1}_{(s,r,z)}$ where the operator $\Psi(s, r, z, D_s, D_r, D_z)$ is microlocally elliptic, the operator

$$\widetilde{\mathcal{M}} := \mathcal{M} - [\mathcal{M}, \Psi(s, r, z, D_s, D_r, D_z)] \Psi(s, r, z, D_s, D_r, D_z)^{-1}$$

is an annihilator of $\widetilde{K}^*\psi d$ to all orders.

**Corollary 8.1.** *A pseudodifferential annihilator of the data is given by*

$$W = K\widetilde{\mathcal{M}}\widetilde{K}^*.$$

Note that $W = W[c_0]$ depends on the background medium. The semi-norm $\|W[c_0]d\|$ can be viewed as the wave equation analog of the differential semblance functional [19].

CHRISTIAAN C. STOLK, THE MATHEMATICAL SCIENCES RESEARCH INSTITUTE, 1000 CENTENNIAL DRIVE #5070, BERKELEY CA 94720-5070 USA
  *E-mail address*: cstolkcaam.rice.edu

MAARTEN V. DE HOOP, THE MATHEMATICAL SCIENCES RESEARCH INSTITUTE, 1000 CENTENNIAL DRIVE #5070, BERKELEY CA 94720-5070 USA
  *E-mail address*: mdehoopMines.EDU